\documentclass[12pt, oneside]{article}   	
\usepackage{geometry}                		
\usepackage{graphicx}	
\usepackage[mathscr]{euscript}			
\usepackage{amssymb}
\def\wt{\widetilde}

 \newtheorem{thm}{Theorem}[section]
 \newtheorem{cor}[thm]{Corollary}
 \newtheorem{lem}[thm]{Lemma}
 \newtheorem{prp}[thm]{Proposition}
 \newtheorem{exm}[thm]{Example}
 \newtheorem{dfn}[thm]{Definition}
 \newtheorem{rmk}[thm]{Remark}
\newenvironment{prf}{\noindent {\it Proof} \ }{\hfill $\Box$}

\title{On stability of subelliptic harmonic maps withpotential}
\author{}

\begin{document}
\title{On stability of subelliptic harmonic maps with potential
\footnotetext{The first author is supported by NSFC Grant No. 12001360. The second author is supported by NSFC Grant No. 11771087, No. 12171091 and LMNS, Fudan. The third author is supported by NSFC Grant No. 12101411.}
}
\author{$\mbox{Tian \ Chong }$ \hspace{12pt} $\mbox{Yuxin\ Dong}$\hspace{12pt} $\mbox{Guilin \ Yang \footnote{Correspondence author}}$ \hspace{12pt}}
\date{}
\maketitle

\begin{abstract}
\textbf{Abstrct:}
In this paper, we investigate the stability problem of subelliptic harmonic maps with potential. First, we derive the first and second variation formulas for subelliptic harmonic maps with potential. As a result, it is proved that a subelliptic harmonic map with potential is stable if the target manifold has nonpositive curvature and the Hessian of the potential is nonpositive definite. We also give Leung type results which involve the instability of subelliptic harmonic maps with potential when the target manifold is a sphere of dimension $\geq 3$.
\end{abstract}

\section{Introduction}

Subellipitc harmonic maps (with potential) from sub-Riemannian manifolds are a natural generalization of harmonic maps (with potential ) from Riemannian manifolds (cf. \cite{Do21}, \cite{DLY22}). Recall that Jost-Xu \cite{JX98} first introduced subelliptic harmonic maps associated with a H\"ormander system of vector fields on a domain of $\mathbb{R}^n$. As a global formulation of Jost-Xu's subelliptic harmonic maps, Barletta et al. \cite{BDU01} introduced subelliptic harmonic maps from strictly pseudoconvex CR manifold into Riemannian manifolds, which are referred to as pseudoharmonic maps in their paper. In particular, they also considered the stability problem for pseudoharmonic maps. In \cite{BD04}, Barletta and Dragomir introduced the notion of pseudoharmonic maps with potential and studied the stability problem for these critical maps from a compact strictly pseudoconvex CR manifold too. In \cite{Do21}, the second author considered a general subelliptic harmonic maps from sub-Riemannian manifolds and established some Eells-Sampson type results. Later, the authrs in \cite{DLY22} generalized the results of \cite{Do21} to the case of subelliptic harmonic maps with potenital.

In this paper, we manily investigate the stability problem for subelliptic harmonic maps with potential. First, we deduce the first and second variation formulas for subelliptic harmonic maps with potential. Using the second variation formula, it is proved that a subelliptic harmonic map with potential is stable if the target manifold has nonpositive curvature and the Hessian of the potential is nonpositive definite. When the target manifold is a sphere of dimension $\geq 3$, we also give some Leung type instability results for subelliptic harmonic map with potential, which generalize related results in \cite{Le82, Ch00} for harmonic maps (with potential).

\section{The first and second variation formulas}

Let $(M^{m+d},g)$ be a compact connected $(m+d)$-dimensional Riemannian manifold and let $H$ be a $m$-dimensional subbundle of $TM$. Set $g_H=g|_{H}$. Then $(M,H,g_H)$ is a sub-Rieamnnian manifold. Let $V$ be the orthogonal complement of $H$ with respect to $g$. Then the tangent bundle $TM$ has the following orthogonal decomposition:
$$
TM=H\oplus V.
$$
The distribution $H$ (resp. $V$) is referred to as the horizontal (resp. vertical) distribution or bundle of $(M,H,g_H; g)$. Set $g_V=g|_{V}$. We can write
$$
g=g_H+g_V.
$$
Let $\nabla$ denote the Levi-Civita connection of $(M,g)$. Then we have the induced connections $\nabla^H$ and $\nabla^V$ on $H$ and $V$ respectively as follows:
\begin{eqnarray*}
\nabla^H_X Y=\pi_H(\nabla_X Y), \quad \nabla^V_X Z=\pi_V(\nabla_X Z),
\end{eqnarray*}
for any $X\in TM$, $Y\in \Gamma(H)$, $Z\in\Gamma(V)$.
Choose a local orthonormal frame field $\{e_A\}_{A=1}^{m+d}$ on an open domain of $(M,g)$ such that $span\{e_i\}_{i=1}^m=H$, and thus $span\{e_{\alpha}\}_{\alpha=m+1}^{m+d}=V$. Such a frame field is referred to as \emph{an adapted frame field} for $(M, H,g_{H};g)$. Let $\{\omega^A\}_{A=1}^{m+d}$ be the dual frame field of $\{e_A\}_{A=1}^{m+d}$.
We shall use the following convention on the ranges of indices:
\begin{eqnarray*}
1\leq A,B,C,\cdots\leq m+d;\\
1\leq i,j,k,\cdots\leq m;\\
m+1\leq \alpha,\beta,\gamma,\cdots\leq m+d.
\end{eqnarray*}

Let $(N^n,h)$ be a Riemannian manifold with the Levi-Civita connection $\widetilde{\nabla}$.
For a smooth map $f:(M^{m+d},H, g_{H};g)\rightarrow (N^n, h)$, we define
\begin{eqnarray*}
df_H=df\circ \pi_H,\\
df_V=df\circ \pi_V,
\end{eqnarray*}
where $\pi_H: TM\rightarrow H$ and $\pi_V:TM\rightarrow V$ are the projection morphisms.
Then we can view $df_H$ and $df_V$ as sections of $T^*M\otimes f^{-1}TN$.
The second fundamental forms for $df_H$ and $df_V$ with respect to $(\nabla, \wt{\nabla}^f)$ are defined respectively by
$$
\beta_H(f)(X,Y)=\widetilde{\nabla}_Y^f (df_H(X))-df_H(\nabla_Y X)
$$
and
$$
\beta_V(f)(X,Y)=\widetilde{\nabla}_Y^f (df_V(X))-df_V(\nabla_Y X)
$$
for any $X,Y\in\Gamma(TM)$, where $\widetilde{\nabla}^f$ denotes the pull-back connection of the Levi-Civita coonection $\widetilde{\nabla}$ on $N$ by $f$.

In terms of the connection $\wt{\nabla}^f$ on $f^{-1}TN$, we may introduce the exterior covariant differential operator $d:\mathcal{A}^p(\xi)\rightarrow \mathcal{A}^{p+1}(\xi)$ and its codifferential operator $d^{*}:\mathcal{A}^p(\xi)\rightarrow \mathcal{A}^{p-1}(\xi)$, where $\xi=f^{-1}TN$ (cf. \cite{EL83}). The Laplacian $\Delta$ is defined on $\xi$-valued differential forms by
\begin{equation}
\Delta=d\circ d^{*}+ d^{*}\circ d:\mathcal{A}^{p}(\xi)\rightarrow \mathcal{A}^{p}(\xi).
\end{equation}

\begin{thm}
(Weitzenb\"ock formula (1.30) of \cite{EL83}) For any $\sigma\in \mathcal{A}^1(\xi)$, we have
\begin{equation}
\Delta \sigma=-trace\wt{\nabla}^2\sigma+S(\sigma),
\end{equation}
where $S(\sigma)(X)=-[\wt{R}(X, e_B)\sigma](e_B)=-\wt{R}(df(X), df(e_B))(\sigma(e_B))+\sigma(R(X, e_B)e_B)$.
\end{thm}

By the Weitzenb\"ock formaula and letting $\sigma=df_H$, we have
\begin{eqnarray}
\nonumber \Delta df_H &=& -trace\wt{\nabla}^2 df_H -\wt{R}(df(\cdot),df(e_B))df_H(e_B)+df_H(Ric^M(\cdot))\\
&=&-trace\wt{\nabla}^2 df_H-\wt{R}(df(\cdot), df(e_i))df_H(e_i)+df_H(Ric^M(\cdot)),
\end{eqnarray}
where $\{e_A\}$ is an adapted frame field for $(M,H,g_H;g)$.

\begin{lem}\label{lem2.2}
For any given $W\in f^{-1}TN$, we have
\begin{equation}
\int_M [e_i\langle W, df(e_i)\rangle-\langle W, df_H(\nabla_{e_i}e_i)\rangle]dv_g=\int_M \langle W, df(\zeta)\rangle dv_g,
\end{equation}
where $\zeta=\pi_H(\nabla_{e_{\alpha}}e_{\alpha})$.
\end{lem}

\begin{prf}
Define a $1$-form $\theta$ on $M$ by
$$
\theta(X)=\langle W, df_H(X)\rangle,
$$
for any $X\in TM$. The codifferential of $\theta$ is given by
\begin{eqnarray*}
d^{*}\theta&=& -(\nabla_{e_A}\theta)(e_A)\\
&=& -[e_A (\theta(e_A))-\theta (\nabla_{e_A}e_A)]\\
&=& -[e_i(\theta(e_i))-\theta(\nabla_{e_i}e_i)]+\theta(\nabla_{e_{\alpha}}e_{\alpha})\\
&=&-[e_i\langle W, df_H(e_i)\rangle-\langle W,df_H(\nabla_{e_i}e_i)]+\langle W, df(\zeta)\rangle,
\end{eqnarray*}
Then the result follows from the divergence theorem.
\end{prf}

Let $f:(M^{m+p},H, g_{H};g)\rightarrow (N^n, h)$ ba a smooth map and $G:N \rightarrow \mathbb{R}$ a smooth function. We can define the horizontal energy with potential as
$$
E(f)=\int_M [e_H(f)-G(f)]dv_g,
$$
where $e_H(f)=\frac{1}{2}|df_{H}|^2$ is the horizontal energy density of $f$, and $dv_g$ is the Riemannian volume element of $g$.

\begin{dfn}
A critical map of $E(f)$ is called a subelliptic harmonic map with potential.
\end{dfn}

\begin{rmk}\label{rem2.4}
If $G$ is a constant function on $N$, then a critical point of $E(f)$ is a subelliptic harmonic map (cf. \cite{Do21}).
\end{rmk}

\begin{thm}\label{thm4}
Let $\{f_t\}_{|t|<\epsilon}$ be a variation of $f:(M,H,g_H; g)\rightarrow (N,h)$ (with $f_0=f$). Suppose the variation vector field $V=\frac{\partial f_t}{\partial t}|_{t=0}\in \Gamma(f^{-1}TN)$ has compact support. Then we have the first variation formula of the horizontal energy with potential
\begin{equation}
\frac{d}{dt}E(f_t)\Big|_{t=0}=-\int_M \langle V, \tau_{HG}(f) \rangle dv_g,
\end{equation}
where $\tau_{HG}(f)=\beta_H(f)(e_i,e_i)-df(\zeta)+\widetilde{\nabla}G(f)$.
\end{thm}

\begin{prf}
Define a map $F:M\times (-\epsilon,\epsilon)\rightarrow N$ by $F(x,t)=f_t(x)$. Clearly
$$
\widetilde{\nabla}_{\frac{\partial }{\partial t}} dF(X)=\widetilde{\nabla}_{X}dF(\frac{\partial }{\partial t}).
$$
Using the above formula, we derive that
\begin{eqnarray}
\frac{d}{dt}E(f_t)&=&\int_M[ \langle\wt{\nabla}_{e_i}dF(\frac{\partial }{\partial t}),dF(e_i)\rangle-\langle \wt{\nabla}G(F(\cdot, t)),\frac{dF(\cdot, t)}{dt} \rangle]dv_g. \label{2.7}
\end{eqnarray}
Putting $t=0$ in the above equation, we have
\begin{eqnarray}
\nonumber \frac{d}{dt}E(f_t)\Big|_{t=0}&=&\int_M [\langle\wt{\nabla}_{e_i}V, df(e_i)\rangle-\langle\wt{\nabla}G(f), V \rangle]dv_g\\
\nonumber &=&\int_M [e_i \langle V, df(e_i)\rangle-\langle V, \wt{\nabla}_{e_i}df(e_i)\rangle -\langle\wt{\nabla}G(f), V \rangle ]dv_g\\
\nonumber &=& \int_M [e_i\langle V, df(e_i)\rangle -\langle V, \beta_H(f)(e_i, e_i)+ df_{H}(\nabla_{e_i}e_i)+\wt{\nabla}G(f)\rangle] dv_g\\
\nonumber &=&-\int_M [\langle V, \beta_H(f)(e_i, e_i)-df(\zeta)+\wt{\nabla}G(f)\rangle] dv_g,
\end{eqnarray}
where the last equality is due to Lemma \ref{lem2.2}.
\end{prf}

Thus $f$ is a subelliptic harmonic map with potential if and only if it satisfies the Euler-Lagrange equation
\begin{equation}
\tau_{HG}(f)=\beta_H(f)(e_i,e_i)-df(\zeta)+\wt{\nabla}G(f)=0.
\end{equation}

\begin{thm}
Let $f, f_t, F, V$ and $E(f_t)$ be as in Theorem \ref{thm4}. Then
\begin{eqnarray}
\nonumber \frac{d^2E(f_t)}{dt^2}|_{t=0}&=&\int_M [\langle \wt{\nabla}_{e_i}V, \wt{\nabla}_{e_i}V \rangle-\langle \wt{R}(df(e_i),V)V, df(e_i)\rangle]dv_g\\
&&-\int_M \big[\langle \xi, \tau_{HG}(f) \rangle-\wt{\nabla}^2G(f)(V,V)\big] dv_g,
\end{eqnarray}
where $\xi= \wt{\nabla}_{\frac{\partial }{\partial t}}dF (\frac{\partial }{\partial t})\Big|_{t=0}\in \Gamma(f^{-1}TN)$.
\end{thm}

\begin{prf}
It follows from (\ref{2.7}) that
\begin{eqnarray*}
\frac{d^2 E(f_t)}{dt^2}&=&\int_M [\langle\wt{\nabla}_{\frac{\partial}{\partial t}}\wt{\nabla}_{e_i}dF(\frac{\partial }{\partial t}),dF(e_i) \rangle+\langle\wt{\nabla}_{e_i}dF(\frac{\partial }{\partial t}),\wt{\nabla}_{\frac{\partial }{\partial t}}dF(e_i) \rangle\\
&& -\langle \wt{\nabla}_{\frac{\partial}{\partial t}}(\wt{\nabla}G(F(\cdot, t))),\frac{dF}{dt} \rangle- \langle \wt{\nabla}G(F(\cdot, t)), \wt{\nabla}_{\frac{\partial}{\partial t}}\frac{dF}{dt} \rangle]dv_g\\
&=&\int_M [\langle -\wt{R}(e_i, \frac{\partial }{\partial t})dF(\frac{\partial }{\partial t})+\wt{\nabla}_{e_i}\wt{\nabla}_{\frac{\partial }{\partial t}}dF(\frac{\partial }{\partial t}), dF(e_i)\rangle+\langle\wt{\nabla}_{e_i}dF(\frac{\partial }{\partial t}),\wt{\nabla}_{e_i}dF(\frac{\partial }{\partial t})\rangle\\
&& - \wt{\nabla}^2G(F(\cdot, t))(\frac{dF}{dt},\frac{dF}{dt})- \langle \wt{\nabla}G(F(\cdot, t)), \wt{\nabla}_{\frac{\partial}{\partial t}}\frac{dF}{dt} \rangle]dv_g.
\end{eqnarray*}
Letting $t=0$ in the above equation, we have
\begin{eqnarray*}
&&\frac{d^2 E(f_t)}{dt^2}\Big|_{t=0}\\
&=& \int_M [\langle -\wt{R}( df(e_i), V)V, df(e_i)\rangle +e_i\langle \wt{\nabla}_{\frac{\partial }{\partial t}}\frac{df_t}{dt}\big|_{t=0}, df(e_i)\rangle
-\langle \wt{\nabla}_{\frac{\partial }{\partial t}}\frac{df_t}{dt}\big|_{t=0}, \wt{\nabla}_{e_i}df(e_i)\rangle\\
&&+\langle\wt{\nabla}_{e_i}V , \wt{\nabla}_{e_i}V\rangle-\widetilde{\nabla}^2G(f)(V,V) -\langle \wt{\nabla}G(f), \wt{\nabla}_{\frac{\partial}{\partial t}}\frac{df_t}{dt}|_{t=0} \rangle]dv_g\\
&=&\int_M [\langle \wt{\nabla}_{e_i}V , \wt{\nabla}_{e_i}V \rangle -\langle \wt{R}( df(e_i), V)V, df(e_i)\rangle-\langle \wt{\nabla}_{\frac{\partial }{\partial t}}\frac{df_t}{dt}\big|_{t=0}, \beta_H(f)(e_i,e_i)\rangle\\
&&+\langle \wt{\nabla}_{\frac{\partial }{\partial t}}\frac{df_t}{dt}\big|_{t=0}, df(\zeta)\rangle-\widetilde{\nabla}^2G(f)(V,V) -\langle \wt{\nabla}G(f), \wt{\nabla}_{\frac{\partial}{\partial t}}\frac{df_t}{dt}|_{t=0} \rangle]dv_g,
\end{eqnarray*}
where the last equality is obtained by applying Lemma 1.1.
\end{prf}

\begin{cor}
Let $f:M\rightarrow N$ be a subelliptic harmonic map with potential. Then
\begin{equation}
\nonumber \frac{d^2E(f_t)}{dt^2}\Big|_{t=0}=\int_M \langle [\wt{\nabla}_{e_i}V, \wt{\nabla}_{e_i}V \rangle-\langle \wt{R}(df(e_i),V)V, df(e_i)\rangle-\wt{\nabla}^2G(f)(V,V)]dv_g.
\end{equation}
\end{cor}

\section{stable subelliptic harmonic maps with potential}

Suppose $f:M\rightarrow N$ is a subelliptic harmonic map with potential. Denote $$
I_{HG}(V,V)=\int_M [\langle \wt{\nabla}_{e_i}V, \wt{\nabla}_{e_i}V \rangle-\langle \wt{R}(df(e_i),V)V, df(e_i)\rangle-\wt{\nabla}^2G(f)(V,V)]dv_g.
$$

\begin{dfn}
A subelliptic harmonic map with potential $f:M\rightarrow N$ is called stable if $I_{HG}(V,V)\geq 0$, for any $V\in \Gamma(f^{-1}TN)$.
\end{dfn}

\begin{prp}
If the target manifold $N$ has nonpositive curvature and the Hessian of $G$ is nonpositive definite (i.e., $\wt{\nabla}^2 G\leq 0$), then $f$ is stable.
\end{prp}

Next, we shall consider the case that the target manifold is a sphere $S^n$.
Let $a$ be a fixed vector in $\mathbb{R}^{n+1}$, define a function on $S^n$ by setting $\varphi(y)=\langle a, y \rangle$ for any $y\in S^n$, where $\langle \cdot, \cdot \rangle $ is also used to denote the inner product in $\mathbb{R}^{n+1}$. Construct a conformal vector field $v=\wt{\nabla}\varphi$ on $S^n$, where $\widetilde{\nabla}$ is the Levi-Civita connection on $S^n$. At any $y\in S^n$, choose a local orthonormal frame fields $\{\tilde{e}_I\}$ such that $(\wt{\nabla}_{\tilde{e}_I}\tilde{e}_J)(y)=0$. It is known that
$$
v=\langle \alpha, \tilde{e}_I\rangle \tilde{e}_I, \quad \wt{\nabla}_X v=-\varphi X,
\quad \wt{\nabla}^2 v=-v,
$$
for any $X\in \Gamma(TS^n)$.
Let $f:M\rightarrow S^n$ be a subelliptic harmonic map with potential. Then we have
\begin{eqnarray}
\nonumber && I_{HG}(v,v)\\
\nonumber &=& \int_M [\langle \wt{\nabla}_{e_i}v, \wt{\nabla}_{e_i}v \rangle
-\langle \wt{R}(df(e_i), v)v, df(e_i) \rangle-\wt{\nabla}^2G(f)(v,v)]dv_g\\
\nonumber &=&\int_M \{2(\varphi^2 (f))e_H(f)-[2e_H(f)\langle v,v \rangle - \langle df(e_i), v \rangle^2] -\wt{\nabla}^2G(f)(v,v)\}dv_g\\
&=& \int_M \{2e_H(f)[\langle \alpha, f(x) \rangle^2-\sum_I\langle \alpha, \tilde{e}_I \rangle^2]+\langle df(e_i), v \rangle^2
-\wt{\nabla}^2G(f)(v,v)\}dv_g.\label{3.1}
\end{eqnarray}
Denote $\{\epsilon_s\}_{s=1}^{n+1}$ the standard orthonormal basis in $\mathbb{R}^{n+1}$. In (\ref{3.1}), we choose $a=\epsilon_s$ and $v_s=\langle \epsilon_s, \tilde{e}_I \rangle \tilde{e}_I$, $s=1,2,\cdots, n+1$, and compute the sum
\begin{eqnarray}
\nonumber && \sum_{s=1}^{n+1}I_{HG}(v_s, v_s)\\
\nonumber &=& \int_M 2e_H(f)\sum_{s}[\langle \epsilon_s, f(x) \rangle^2-\sum_{I}\langle \epsilon_s, \tilde{e}_I \rangle^2]+\sum_{i,s}\langle df(e_i), \epsilon_s \rangle^2-\sum_s \wt{\nabla}^2G(f)(v_s,v_s)\\
&=& 2(2-n)\int_M e_H(f)-\sum_s \int_M\wt{\nabla}^2G(f)(v_s,v_s).
\end{eqnarray}

\begin{thm}\label{thm3.3}
Suppose $f:M\rightarrow S^n$ is a nonconstant subelliptic harmonic map with potential from a compact sub-Riemannian manifold $M$ into $S^n$ ($n\geq 3$) with $f(M)\subset U$, where $U$ is a domain of $S^n$. If $G:U\rightarrow \mathbb{R}$ is a convex function, i.e., $\wt{\nabla}^2G(V,V)\geq 0$ for any vector field $V\in \Gamma(TU)$, then $f$ is unstable.
\end{thm}

\begin{exm}
If $U$ in theorem \ref{thm3.3} is the upper half sphere $S^n_{+}$, the potential function $G$ can be selected as follows.\\
(1) By the Hessian comparison theorem, it is easy to find that $r^2(y)$ in convex on $S^n_{+}$, where $r(y)=dist(y, p)$ is the Riemannian distance function on $S^n$ relative to the north point $p$.\\
(2) If the potential function is taken by $G(y)=-\langle \epsilon_{n+1}, y\rangle$ on $S^n$, then we have $\wt{\nabla}^2 G(f)(V,V)=-|V|^2\cdot G(f)$. And hence $G(y)$ is convex on $S^n_{+}$.
\end{exm}

\begin{cor}
If $n\geq 3$, then any nonconstant subelliptic harmonic map with constant potential function $G$ from a compact sub-Riemannian manifold $M$ into $S^{n}$ is unstable.
\end{cor}

From Remark \ref{rem2.4}, we know that there is no nonconstant subelliptic harmonic map from compact sub-Riemannian manifolds into $n$-dimensional  spheres with $n\geq 3$.

\begin{flushleft}
Tian Chong\\
Department of Mathematics\\
Shanghai Polytechnic University, Shanghai 201209, P.R. China\\
E-mail address: chongtian@sspu.edu.cn
\end{flushleft}

\begin{flushleft}
Yuxin Dong\\
School of Mathematical Science\\
and\\
Laboratory of Mathematics for Nonlinear Science\\
Fudan University, Shanghai 200433, P.R. China\\
E-mail address: yxdong@fudan.edu.cn
\end{flushleft}

\begin{flushleft}
Guilin Yang\\
School of Statistics and Mathematics\\
Shanghai Lixin University of Accounting and Finance, Shanghai 201209 P.R. China\\
E-mail address: glyang@lixin.edu.cn
\end{flushleft}

\end{document}